\documentclass[preprint,onecolumn,number]{elsarticle}
\usepackage{amsmath,amsthm,amssymb}
\usepackage[mathscr]{euscript}
\usepackage[ruled,vlined]{algorithm2e}
\usepackage{hyperref}

\newtheorem{thm}{Theorem}
\newtheorem{defn}{Definition}
\newtheorem{prop}{Proposition}
\newtheorem{rem}{Remark}
\newtheorem{ex}{Example}

\newcommand{\lsm}[1]{\mbox{\sf{l-SJ}}(#1)}
\newcommand{\rsm}[1]{\mbox{\sf{r-SJ}}(#1)}
\newcommand{\uu}[1]{{\mathcal U}_{#1}[\delta,\ddt]}

\newcommand{\uud}[1]{{\mathcal U}_{#1}\left( \delta\right)[\ddt]}
\newcommand{\mm}[2]{{\mathcal M}_{#1,#2}[\delta,\ddt]}
\newcommand{\mmd}[2]{{\mathcal M}_{#1,#2}\left( \delta\right) [\ddt]}
\newcommand{\mmdd}[2]{{\mathcal M}_{#1,#2}\left( \delta_{1},\delta_{2}\right) [\ddt]}
\newcommand{\mmp}[2]{{\mathcal M}_{#1,#2}[\delta,\pi^{-1},\ddt]}

\newcommand{\wt}[1]{\widetilde{#1}}

\newcommand{\ddt}{{\frac{d}{dt}}}
\newcommand{\kk}{{\mathfrak K}}
\newcommand{\NN}{{\mathbb N}}
\newcommand{\RR}{{\mathbb R}}
\newcommand{\ds}{\displaystyle}
\newcommand{\diag}[1] {\mathsf{diag}\{ #1 \}}
\newcommand{\lcm}[1]{\mbox{\sf{LCM}}(#1)}
\begin{document}
\begin{frontmatter}
\title{On the Computation of $\pi$-Flat Outputs for\\Linear Time-Delay Systems\tnoteref{label1}}
\tnotetext[label1]{This work has been partly done while the authors were participating in the Bernoulli Program ``Advances in the Theory of Control, Signals, and Systems, with Physical Modeling'' of the Bernoulli Center, EPFL, Switzerland.}

\author{F.~Cazaurang}
\ead{franck.cazaurang@ims-bordeaux.fr}
\address{University of Bordeaux, IMS Laboratory, CNRS, UMR 5218, 351 cours de la Lib\'{e}ration, 33405 Bordeaux, France.}

\author{J.~L\'{e}vine}
\ead{jean.levine@mines-paristech.fr}
\address{Centre Automatique et Syst\`{e}mes (CAS), Unit\'{e} Math\'{e}matiques et Syst\`{e}mes, MINES-ParisTech, 35 rue Saint-Honor\'{e}, 77300 Fontainebleau, France.}

\author{V.~Morio\tnoteref{label2}}
\tnotetext[label2]{Now with DGA/DT/LRBA/SDT/MAN/PSP, BP 914, 27207 Vernon Cedex, France.}
\ead{vincent.morio@dga.defense.gouv.fr}
\address{University of Bordeaux, IMS Laboratory, CNRS, UMR 5218, 351 cours de la Lib\'{e}ration, 33405 Bordeaux, France.}

\begin{abstract}
This paper deals with the concepts of $\pi$-flatness and $\pi$-flat output for linear time-varying delay systems. These notions, introduced and developed by several authors during the last decade, may be, roughly speaking, defined as follows: a $\pi$-flat system is a system for which all its variables may be expressed as functions of a particular output $y$, a finite number of its successive time derivatives, time delays, and predictions, the latter resulting from the advance operator $\pi^{-1}$, $\pi$ being a polynomial of $\delta$, the delay operator. 
Thanks to standard polynomial algebraic tools, and in particular the Smith-Jacobson decomposition of polynomial matrices, we obtain a simple and easily computable characterization of $\pi$-flatness in terms of \emph{hyper-regularity} of the system matrices and deduce a constructive algorithm for the computation of $\pi$-flat outputs. Some examples are provided to illustrate the proposed methodology.
\end{abstract}

\end{frontmatter}

\section{Introduction}
Differential flatness, roughly speaking, means that all the variables of an under-determined system of differential equations can be expressed as functions of a particular output, called flat output, and a finite number of its successive time derivatives (\cite{Martin_92,Fliess_95,Fliess_99}, see also  \cite{Ramirez_04,Levine_09,Levine_11} and the references therein). 

For time-delay systems and more general classes of infinite-dimensional systems, extensions of this concept have been proposed  and thoroughly discussed in \cite{Mounier_95,Fliess_96,Petit_00,Rudolph_03}. In a linear context, relations with the notion of system parameterization \cite{Pommaret_99,Pommaret_01} and, in the behavioral approach of \cite{Polderman_98}, with \emph{latent variables} of \emph{observable image representations} \cite{Trentelman_04}, have been established. Other theoretic approaches have been proposed e.g. in \cite{Rocha_97,Chyzak_05}. Interesting control applications of linear time-delay systems may be found in \cite{Mounier_95,Petit_00,Rudolph_03}.

Characterizing differential flatness and flat outputs has been an active topic since the beginning of this theory. The interested reader may find a historical perspective of this question in \cite{Levine_09,Levine_11}. Constructive algorithms, relying on standard computer algebra environments, may be found e.g. in \cite{Antritter_08} for nonlinear finite-dimensional systems, or \cite{Chyzak_04} for linear systems over Ore algebras. 

The results and algorithm proposed in this paper for the characterization and computation of $\pi$-flat outputs for linear time-delay systems are strongly related to the algebraic framework developed in \cite{Mounier_95,Petit_ecc97,Rudolph_03}. 
More precisely, we study linear time-delay differential control systems, i.e. linear systems of the form $Ax=Bu$, with $x \in \RR^{n}$ the pseudo-state, and $u\in \RR^{m}$ the control, for given integers $m\leq n$, where the entries  of the matrices $A$ and $B$ belong to the ring $\kk[\delta,\ddt]$ of multivariate polynomials of $\delta$, the delay operator, and $\ddt$, the time derivative operator, over the ground field  $\kk$ of meromorphic functions of the variable $t$.
We say that the system $Ax=Bu$ is $\pi$-flat if, and only if, the module generated by the components of $x$ and $u$ over $\kk[\delta,\ddt]$ and satisfying the relations $Ax=Bu$, localized at the powers of a polynomial $\pi\in \kk[\delta]$, is free, and a $\pi$-flat output is a basis of this free module (see \cite{Mounier_95}).

To characterize and compute $\pi$-flat outputs, we propose a methodology based on standard polynomial algebra, generalizing the one used in \cite{Levine_03} for ordinary linear differential systems, by extending the original ring $\kk[\delta,\ddt]$ to the principal ideal ring $\kk(\delta)[\ddt]$ of polynomials of $\ddt$ over the fraction field $\kk(\delta)$, namely the field generated by fractions of polynomials of $\delta$ with coefficients in $\kk$, and finally localize the results of our computations at the powers of a suitable polynomial $\pi$ of $\kk[\delta]$. This approach allows us to use the well-known Smith-Jacobson (or diagonal) decomposition (\cite{Cohn_85,Jacobson_78}) of matrices with entries in the larger ring $\kk(\delta)[\ddt]$ as the main tool to obtain the searched $\pi$-flat outputs.

Following \cite{Levine_03}, in order to work with a smaller set of equations and variables, we eliminate the input variables, leading to an implicit system representation, as opposed to previous approaches (see e.g. \cite{Mounier_95,Petit_00,Rudolph_03,Chyzak_05,Chyzak_07}).
Let us also insist on the fact that the time-varying dependence of the systems under consideration is in the class of meromorphic functions, whereas in \cite{Chyzak_05,Chyzak_07}, this dependence is polynomial with respect to time in order to apply effective Gr\"{o}bner bases techniques. 

The main contributions of this paper are (1) the characterization of $\pi$-flatness in terms of the \emph{hyper-regularity} of the system matrices, (2) yielding an elementary algorithm to compute $\pi$-flat outputs, based on the Smith-Jacobson decomposition of the former matrices. 
In addition, the evaluation of our $\pi$-flatness criterion only relies on computations over the larger ring $\kk(\delta)[\ddt]$.

The paper is organized as follows. The $\pi$-flat output computation problem is described in section~\ref{statement-sec}, as well as the algebraic framework. Then, the main result of the paper is presented in section~\ref{main-sec}. Finally, the proposed methodology is illustrated by some examples  in section~\ref{ex-sec}, and its generalization to multiple delays is outlined on an example of vibrating string, first solved in \cite{Mounier_98}.

\section{Problem Statement}\label{statement-sec}

We consider a linear system governed by the set of time-delay differential equations:
\begin{equation}
A\left(\delta,\ddt \right)x=B\left(\delta,\ddt \right)u,\label{sys_lin_2}
\end{equation}
where $x\in\RR^{n}$ is the pseudo-state, $u\in\RR^{m}$ the input vector, $A$ (resp. $B$) a $n\times n$ (resp. $n\times m$) matrix, whose coefficients are multivariate polynomials of $\delta$ and $\ddt$, with $\ddt$ the differentiation operator with respect to time and $\delta$ the time-delay operator defined by:
\begin{equation}
\delta:f(t)\mapsto \delta f(t)=f(t-\tau),\quad \forall t \in \RR\label{delay_eq}
\end{equation}
where $\tau \in \RR^{+}$ is the delay.

In order to precise the nature of the coefficients $a_{i,j}(\delta,\ddt)$, $i,j=1,\ldots,n$, and $b_{i,j}(\delta,\ddt)$, $i=1,\ldots,n$, $j=1,\ldots,m$, of the matrices $A(\delta,\ddt)$ and $B(\delta,\ddt)$ respectively, some algebraic recalls are needed.

\subsection{Algebraic Framework}\label{alg_sec}
Since we deal with smooth functions of time, a natural field is the \emph{differential field of meromorphic functions on the real line} $\RR$. We call this field the \emph{ground field} and we denote it by $\kk$. 
The previously introduced operators $\delta$ and $\ddt$ satisfy the following rules:
\begin{equation}
\ddt \left( \alpha(t)\cdot \right)=\alpha(t)\ddt+\dot{\alpha}(t)\cdot,\quad\delta \left(\alpha(t)\cdot\right)=\alpha(t-\tau)\delta,\quad\ddt\delta=\delta\ddt\nonumber
\end{equation}
for every time function $\alpha$ belonging to $\kk$. The set of multivariate polynomials of these operators, namely polynomials of the form 
$$\sum_{k,l~\rm{finite}}\alpha_{k,l}(t)\frac{d^{k}}{dt^{k}}\delta^{l},\quad \alpha_{k,l}\in\kk$$ 
is a \emph{skew commutative ring} \cite{McConnell_00,Shafarevich_97}, denoted by $\kk[\delta,\ddt]$.
The coefficients $a_{i,j}$ (resp. $b_{i,j}$) of the matrix $A$ (resp. $B$)  of system (\ref{sys_lin_2}) are supposed to belong to $\kk[\delta,\ddt]$, thus making system (\ref{sys_lin_2}) a linear time-varying time-delay differential system, whose coefficients are meromorphic functions with respect to time.

\subsubsection{System Module, Freeness}\label{module-sec}
To system (\ref{sys_lin_2}) is associated the so-called \emph{system module}, noted $\Lambda$. More precisely, following \cite{Fl-scl, Mounier_95}, let us consider a non zero, but otherwise arbitrary,  pair $(\xi,\nu)= (\xi_{1},\ldots, \xi_{n},\nu_{1},\ldots, \nu_{m})$ and the free module\footnote{For more details on rings and modules, the reader may refer to \cite{Cohn_85}.}, denoted by $[\xi,\nu]$, generated by all possible linear combinations of $\xi$ and $\nu$ with coefficients in $\kk[\delta,\ddt]$. Next, we set $\theta=A\xi-B\nu$ and construct the submodule $[\theta]$ of $[\xi,\nu]$ generated by the components of the vector $\theta$.
The system module $\Lambda$ is, by definition, the quotient module $\Lambda =[\xi,\nu]/[\theta]$.

In \cite{Mounier_95}, in the context of commutative polynomial rings, the notion of projective (resp. torsion-free) \emph{controllability} of a time-invariant system, i.e. a system of the form (\ref{sys_lin_2}) with ground field $\kk = \RR$, is defined as the projective (resp. torsion) freeness of $\Lambda$, and shown to generalize the well-known Kalman controllability criterion to linear time-invariant differential delay systems. Moreover, as a consequence of a theorem of Quillen and Suslin, solving a conjecture of Serre (see e.g. \cite{Eisenbud_94,Lam_78}), $\Lambda$ is free if and only if it is projective free.
If $F$ is a finite-dimensional presentation matrix of $\Lambda$, the latter module $\Lambda$ is projective free if $F$ is right-invertible, i.e. there exists a matrix $T$ over $\kk[\delta,\ddt]$ such that $FT = I$. 

This approach has been generalized to modules over the Weyl algebras by Quadrat and Robertz \cite{Quadrat_07}, based on a theorem of Stafford \cite{Stafford_78}, (see algorithmic versions of this result in \cite{Hillebrand_01,Leykin_04}).

In both time-invariant and time-varying cases, systems whose module is free are called \emph{flat} (\cite{Mounier_95,Petit_ecc97,Rudolph_03}). Nevertheless, only few systems have a free system module, thus motivating the weaker notion of $\pi$-flatness: we say that the system is $\pi$-flat, or that its associated module is $\pi$-free (\cite{Mounier_95,Rudolph_03}), if, and only if, there exists a polynomial $\pi\in \kk[\delta]$, called \emph{liberation polynomial} (\cite{Mounier_95,Rudolph_03}), such that the module $\kk[\delta, \pi^{-1},\ddt] \otimes_{\kk[\delta, \ddt]} \Lambda$, i.e. the set of elements of the form $\sum_{i\in I} \pi^{-i} a_{i}\xi_{i}$ with  $I$ arbitrary subset of $\NN$, $a_{i}\in \kk[\delta, \ddt]$ and $\xi_{i} \in \Lambda$ for all $i\in I$, called the \emph{system module localized at the powers of $\pi$}, is free.
In other words, $\pi$-flatness means that the state and input can be expressed in terms of the $\pi$-flat output, a finite number of its time derivatives and delays, and advances corresponding to powers of the inverse  operator $\pi^{-1}$.

In the sequel, we will also use the extension, as announced, of the ground field $\kk$ to $\kk(\delta)$, the fraction field generated by rational functions of $\delta$ with coefficients in $\kk$. The system module over this field extension is  $\kk(\delta)[\ddt] \otimes_{\kk[\delta, \ddt]} \Lambda$. Indeed, freeness (in any sense) of the latter module does not imply freeness (in any sense) of the original system module $\Lambda$ (see e.g. \cite{Mounier_95}).


\subsubsection{Polynomial Matrices, Smith-Jacobson Decomposition, Hyper-regularity}\label{SJdecomp-sec}
The matrices of size $p\times q$ whose entries are in $\kk[\delta,\ddt]$ generate a \emph{module} denoted by $\mm{p}{q}$. The matrix $M\in\mm{n}{n}$ is said to be \emph{invertible}\footnote{Note that the $\kk[\delta,\ddt]$-independance of the $n$ columns and rows of $M$ is not sufficient for its invertibility. Its inverse, denoted by $N$, has to be polynomial too.} if there exists a matrix $N\in\mm{n}{n}$ such that $MN=NM=I_{n}$, where $I_{n}$ is the identity matrix of order $n$. The subgroup of $\mm{n}{n}$ of invertible matrices is called the group of \emph{unimodular matrices} of size $n$ and is noted $\uu{n}$\footnote{It is also often denoted by $GL_{n}(\kk[\delta, \ddt])$}.

Let us give an example of a system of the form (\ref{sys_lin_2}), that will serve as a guideline all along this section to illustrate the various concepts.
\begin{ex}\label{Example_1}
\begin{equation}
Ax \triangleq \left(\begin{array}{cc}
\ddt&-k(t)(\delta-\delta^{2})\\0&\ddt
\end{array}\right)x=
\left(\begin{array}{c}
0\\\delta
\end{array}\right)u \triangleq Bu\label{example_1}
\end{equation}
where $x=(x_{1},x_{2})^{T}$, $u$ is scalar and $k(t)$ a meromorphic function. In other words (\ref{example_1}) reads:
\begin{equation}
\left\{\begin{array}{l}
\dot{x}_{1}(t)=k(t)(x_{2}(t-\tau)-x_{2}(t-2\tau))\\\dot{x}_{2}(t)=u(t-\tau)
\end{array}\right.\label{example_1_bis}
\end{equation}
The coefficients $\ddt$ and $-k(t)(\delta-\delta^{2})=-k(t)\delta (1-\delta)$ are elements of $\kk[\delta,\ddt]$ and the corresponding matrices $A$ and $B$ belong to $\mm{2}{2}$ and $\mm{2}{1}$ respectively. 
\end{ex}
Note that it may be necessary to extend the polynomial ring as shown by the following computation on the previous example: 

Let us express $u$ of (\ref{example_1}) as a function of $x_{1}$. It is straightforward to see that $x_{2}=\ds (1-\delta)^{-1}\delta^{-1}\frac{1}{k}\ddt x_{1}$ and, since $u=\delta^{-1}\dot{x}_{2}$, we immediately get \mbox{$u = \delta^{-2}(1-\delta)^{-1}\left(\ds- \frac{\dot{k}}{k^{2}}+\frac{1}{k}\ddt\right)\ds\ddt x_{1}$}. Denoting by $\pi\triangleq 
(1-\delta)\delta^{2}\in \kk[\delta]$, 
the polynomial $\pi^{-1}\left(\ds - \frac{\dot{k}}{k^{2}}+\frac{1}{k}\ddt\right)\ds\ddt$ lives in $\kk[\delta,\pi^{-1},\ddt]$ the ring of polynomials of $\delta,\pi^{-1}$ and $\ddt$ with coefficients in $\kk$, but not in $\kk[\delta,\ddt]$. This is why we may also introduce matrices over $\kk[\delta,\pi^{-1},\ddt]$ for some given $\pi$ in $\kk[\delta]$. The corresponding module of matrices of size $p\times q$ will be denoted by $\mmp{p}{q}$. More precisely, a matrix $M$ belongs to $\mmp{p}{q}$ if and only if there exists a finite $s\in \NN$ such that $\pi^{s} \cdot M \in \mm{p}{q}$.
\begin{rem}\label{finitesupp-rem}
It may be argued that the previous expression of $u$ in function of $x_{1}$ is not feasible since $(1-\delta)^{-1}=\sum_{j=0}^{+\infty}\delta^{j}$ implies 
$$u(t)=\ds\sum_{j=-2}^{+\infty} \left(- \frac{\dot{k}(t-j\tau)}{k^{2}(t-j\tau)}\dot{x}_{1}(t-j\tau)+\frac{1}{k(t-j\tau)}\ddot{x}_{1}(t-j\tau)\right)$$ 
which involves an infinite number of delayed terms. However, if we deal with motion planning, if $x_{1}$ is chosen constant outside the interval $[t_{0}, t_{1}]$, for some $t_{0}, t_{1} \in \RR$, $t_{0} < t_{1} $, then $\dot{x}_{1}$ and $\ddot{x}_{1}$ are identically equal to zero before $t_{0}$ and after $t_{1}$, and the above series has at most $\lceil \frac{t_{1}-t_{0}}{\tau}\rceil+2$ non vanishing terms, where the notation $\lceil r \rceil$ stands for the least integer upper bound of an arbitrary non negative real $r$. 
\end{rem}  
Unfortunately, $\kk[\delta,\ddt]$ and $\kk[\delta,\pi^{-1},\ddt]$ are not Principal Ideal Domains (see e.g. \cite{McConnell_00,Shafarevich_97,Jacobson_78,Cohn_85}), a property which is essential for our purpose (see the Smith-Jacobson decomposition in Appendix~\ref{SJalgo-sec}). However, if we extend the ground field $\kk$ to the fraction field $\kk(\delta)$, $\kk(\delta)[\ddt]$ is a principal ideal ring of polynomials of $\ddt$. We then construct the modules $\mmd{p}{q}$ of matrices of size $p\times q$ and $\uud{p}$ of unimodular matrices of size $p\times p$ respectively, over $\kk(\delta)[\ddt]$. Note that $\kk(\delta)[\ddt]$ strictly contains $\kk[\delta,\ddt]$ and $\kk[\delta,\pi^{-1},\ddt]$ for every $\pi\in\kk[\delta]$. Therefore, to interpret the results of computations in $\kk(\delta)[\ddt]$, which turn out to be quite simple, and to decide if they belong to a suitable $\kk[\delta,\pi^{-1},\ddt]$, following \cite{Mounier_95}, we have recourse to the notion of \emph{localization} introduced in subsection~\ref{module-sec}. This aspect will be discussed in section~\ref{main-sec}.

Since $\kk(\delta)[\ddt]$ is a Principal Ideal Domain, the matrices of $\mmd{p}{q}$ enjoy the essential property of admitting a so-called \emph{Smith-Jacobson decomposition}\footnote{we adopt here the names of Smith and Jacobson for the diagonal decomposition to remind that it is credited to Smith \cite{Gantmacher_66,Kailath_79} in the commutative context and Jacobson \cite{Jacobson_78,Chyzak_05} for general principal ideal domains.}, or \emph{diagonal decomposition}:
\begin{thm}[Smith-Jacobson decomposition \cite{Cohn_85,Jacobson_78}]\label{SJdecomp-thm}
Let $M\in\mmd{p}{q}$ be an arbitrary polynomial matrix of size $p\times q$. There exist unimodular matrices $U\in\uud{p}$ and $V\in\uud{q}$ such that:
\begin{equation} \label{smith-Jac-decomp} 
UMV=\left\{ 
\begin{array}{ll}
\ds (\Delta_{p}|0_{p,q-p})& \mbox{\textrm{if~~}} p\leq q\vspace{0.5em}\\
\left(\begin{array}{c}\Delta_{q}\\0_{p-q,q}\end{array}\right)&\mbox{\textrm{if~~}} p>q.
\end{array}
\right.
\end{equation}
In both cases, $\Delta_{\sigma}\in\mmd{\sigma}{\sigma}$, $\sigma = p$ or $q$, is a diagonal matrix whose diagonal elements $(d_{1},\ldots,d_{s},0,\ldots,0)$, with $s\leq \sigma$, are such that $d_{i}$ is a nonzero $\ddt$-polynomial for $i=1,\ldots,s$, with coefficients in $\kk(\delta)$, and is a divisor of $d_{j}$ for all $1\leq j\leq i$.
\end{thm}
A constructive algorithm to compute this decomposition may be found in Section~\ref{SJalgo-sec} of the Appendix.

Given an arbitrary matrix $M$, we call $\lsm{M}$ (resp. $\rsm{M}$), the left (resp. right) Smith-Jacobson subset of unimodular matrices $U\in\uud{p}$ (resp. $V\in\uud{q}$) such that there exists $V\in\uud{q}$ (resp. $U\in\uud{p}$) satisfying the decomposition (\ref{smith-Jac-decomp}).

\begin{ex}\label{example_smith}
Consider again system (\ref{example_1}). The Smith-Jacobson decomposition of $B$ is straightforward:
\begin{equation}
U=\left(\begin{array}{cc}0&1\\1&0\end{array}\right),\quad V=1,\quad\Delta_{1}=\delta, \quad UBV= \left(\begin{array}{c}\delta\\0\end{array}\right)\nonumber
\end{equation}
If we want to eliminate $u$ in (\ref{example_1}), using the previous Smith-Jacobson decomposition  of $B$, we first remark that the second line of $U$, which will be denoted by $U_{2}=\left( 1\quad 0\right)$, corresponds to the left projection operator on the kernel of $B$, i.e. $U_{2}B=0$. It suffices then to left multiply $A$ by $U_{2}$ to obtain the implicit form 
\begin{equation}
F(\delta,\ddt)x  \triangleq U_{2}Ax=\left(\begin{array}{cc}\ddt&-k\delta(1-\delta)\end{array}\right)\left(\begin{array}{c}x_{1}\\x_{2}\end{array}\right)=U_{2}Bu = 0\label{implicit_form}
\end{equation}

We may also compute a Smith-Jacobson decomposition of $F$: we first right multiply $F$ by $\left(\begin{array}{cc}0&1\\1&0\end{array}\right)$ to shift the 0-th order term in $\ddt$ to the left, yielding
$$F\left(\begin{array}{cc}0&1\\1&0\end{array}\right)=\left(\begin{array}{cc}-k\delta(1-\delta)&\ddt\end{array}\right)$$ and then, right multiplying the result by 
$$\left(\begin{array}{cc}-\delta^{-1}(1-\delta)^{-1}\frac{1}{k}&\delta^{-1}(1-\delta)^{-1}\frac{1}{k}\ddt\\0&1\end{array}\right)$$ 
leads to $(1\quad 0)$. The Smith-Jacobson decomposition of $F$ is therefore given by
\begin{equation}
U_{F}FV_{F}=(1\quad 0).\label{F_decomp}
\end{equation}
with
\begin{equation}
\begin{aligned} 
V_{F}&=\left(\begin{array}{cc}0&1\\1&0\end{array}\right)\left(\begin{array}{cc}- \delta^{-1}(1-\delta)^{-1}\frac{1}{k}&\delta^{-1}(1-\delta)^{-1}\frac{1}{k}\ddt\\0&1\end{array}\right)\\
&=\left(\begin{array}{cc}0&1\\-\delta^{-1}(1-\delta)^{-1}\frac{1}{k}&\delta^{-1}(1-\delta)^{-1}\frac{1}{k}\ddt
\end{array}\right)
\label{Smith_VF}
\end{aligned}
\end{equation}
and $U_{F}=1$. 

As previously discussed, the Smith-Jacobson decomposition has been computed over the ring $\kk(\delta)\left[\ddt\right]$ but, according to (\ref{Smith_VF}), its result may be expressed in the ring $\kk[\delta, \pi^{-1}_{F},\ddt]$ with $\pi_{F}=(1-\delta)\delta$.

It is also easy to verify that, since the matrix $F$ is a presentation matrix of the system module over $\kk[\delta, \pi^{-1}_{F},\ddt]$, the latter module, according to (\ref{F_decomp}), is isomorphic to any free finitely generated module admitting the matrix $(1\quad 0)$ as presentation matrix, which implies that the localized system module is free (see Section~\ref{module-sec}). Note also that the polynomial $\pi_{F}$ admits non zero roots: every $\tau$-periodic non zero meromorphic function $f$ of the variable $t$ satisfies $\pi_{F}f(t)=f(t-\tau)-f(t-2\tau)=0$. Therefore, the system module is not $\kk[\delta,\ddt]$-free.
\end{ex}

We now introduce a remarkable class of matrices of $\mmd{p}{q}$ called \emph{hyper-regular}. 
\begin{defn}[Hyper-regularity \cite{Levine_09}]
Given a matrix $M\in\mmd{p}{q}$, we say that $M$ is $\kk(\delta)\left[\ddt\right]$-\emph{hyper-regular}, or simply \emph{hyper-regular} if the context is non ambiguous, if and only if, in (\ref{smith-Jac-decomp}), $\Delta_p=I_p$ if $p\leq q$ (resp. $\Delta_q=I_q$ if $p> q$).
\end{defn}

\begin{rem} It is not difficult to prove that a finitely generated module $\Lambda$ over the ring $\kk(\delta)[\ddt]$ whose presentation matrix $F$ is hyper-regular cannot have torsion elements and is therefore free, since $\kk(\delta)[\ddt]$ is a Principal Ideal Ring. However, the system module $\Lambda$, over $\kk[\delta,\ddt]$, does not need to be free, as shown in the previous example. Nevertheless, it will be seen later that there exists a \emph{liberation polynomial} $\pi$, deduced from the Smith-Jacobson decomposition of $F$, such that the $\kk[\delta,\pi^{-1},\ddt]$-system module, associated to the same presentation matrix $F$, is free.
\end{rem}

\begin{ex}\label{example_smith_bis}
Going back to the decomposition of Example~\ref{example_smith}, a variant of this decomposition may be obtained as:
\begin{equation}
U'=\left(\begin{array}{cc}0&\delta^{-1}\\1&0\end{array}\right),\quad V=1,\quad \Delta'_{1}=I_{1}=1\label{Smith_B}
\end{equation}
which proves that $B$ is $\kk(\delta)\left[\ddt\right]$-hyper-regular.

It is also immediate from (\ref{F_decomp}) that $F$ is $\kk(\delta)\left[\ddt\right]$-hyper-regular, or more precisely, $\kk\left[\delta, \pi_{F}^{-1},\ddt\right]$-hyper-regular with $\pi_{F}=(1-\delta)\delta$. Note, on the contrary, that $A$ is not $\kk(\delta)\left[\ddt\right]$-hyper-regular, though $F=U_{2}A$ is. It is easily verified that
$$U_{A}AV_{A}= \Delta_{A}$$
with 
$$\begin{array}{lcl}
U_{A}&=&\left( \begin{array}{cc}
1&0\\
-\frac{\dot{k} }{k}+ \ddt&k\delta(1-\delta)
\end{array}\right) \vspace{0.5em}
\\ 
V_{A}&=&\left( \begin{array}{cc}
0&1\\
-(1-\delta)^{-1}\delta^{-1}\frac{1}{k}&(1-\delta)^{-1}\delta^{-1}\frac{1}{k}\ddt
\end{array}\right)
\end{array}
$$
and 
$$ \Delta_{A}=\left( \begin{array}{cc}
1&0\\
0&(-\frac{\dot{k} }{k} +\ddt )\ddt
\end{array}\right).
$$
Thus, since the second degree $\ddt$-polynomial $(-\frac{\dot{k} }{k} +\ddt)\ddt$ of the diagonal of $\Delta_{A}$ cannot be reduced, $A$ is not hyper-regular.
\end{ex}


%
\subsubsection{Implicit system representation}

One of the applications of the Smith-Jacobson decomposition concerns the possibility of expressing the system (\ref{sys_lin_2}) in implicit form by eliminating the input $u$, which may be useful to work with a smaller number of variables. 

For simplicity's sake, we rewrite system (\ref{sys_lin_2}) $Ax=Bu$.

\begin{prop}\label{impl-prop}
System (\ref{sys_lin_2}) is equivalent to 
\begin{equation}\label{semi_impl_rep}
Fx=0, \quad \Delta_{B}N^{-1}u=(I_{m},0_{m,n-m})MAx
\end{equation}
with 
\begin{equation}\label{Fdef}
F= (0_{n-m,m},I_{n-m})MA,
\end{equation} 
$M\in \lsm{B}$ and $N$ such that 
\begin{equation}\label{Bdecomp}
MBN=\left(\begin{array}{c}\Delta_{B}\\0_{n-m,m}\end{array}\right).
\end{equation}

Moreover, if $B$ is $\kk(\delta)\left[\ddt\right]$-hyper-regular, the explicit form (\ref{sys_lin_2}) admits the implicit representation
\begin{equation}\label{impl_rep}
Fx=0
\end{equation}
with $F$ given by (\ref{Fdef}), and with $\Delta_{B}=I_{m}$ in (\ref{Bdecomp}). In this case, $u$ is deduced from $x$ by 
\begin{equation}\label{udef}
u=N(I_{m},0_{m,n-m})MAx.
\end{equation}
\end{prop}

\begin{proof}
Consider a pair of matrices $M$ and $N$ obtained from the Smith-Jacobson decomposition of $B$, i.e. satisfying (\ref{Bdecomp}). Thus, left-multiplying both sides of system (\ref{sys_lin_2}) by $(0_{n-m,m},I_{n-m})M$, according to (\ref{Fdef}) we get $Fx=0$. On the other hand, multiplying both sides of (\ref{sys_lin_2}) by $(I_{m},0_{m,n-m})M$ we get 
$$(I_{m},0_{m,n-m})MAx= \Delta_{B}N^{-1}u,$$ 
hence the representation (\ref{semi_impl_rep}).

Conversely, if $x$ and $u$ are given by (\ref{semi_impl_rep}), we have
$$MAx=\left( \begin{array}{c} (I_{m},0_{m,n-m})MA\\(0_{n-m,m},I_{n-m})MA\end{array}\right)x = \left( 
\begin{array}{c} \Delta_{B}N^{-1}u\\0 \end{array}\right)= MBu$$
the last equality being a consequence of (\ref{Bdecomp}). Thus, since $M$ is unimodular, the pair $(x,u)$ satisfies $Ax=Bu$, which proves the equivalence.

Finally, if $B$ is $\kk(\delta)\left[\ddt\right]$-hyper-regular, one can replace $\Delta_{B}$ by $I_{m}$ and the second equation of (\ref{semi_impl_rep}) becomes (\ref{udef}). Thus, $u$ is a $\kk\left(\delta\right)[\ddt]$-combination of the components of $x$ and can be eliminated. Therefore, the remaining part (\ref{impl_rep}) is the desired implicit representation of (\ref{sys_lin_2}). The proposition is proven.
\end{proof}

In the sequel, if $B$ is hyper-regular,  we refer to (\ref{impl_rep}) as the implicit representation of system (\ref{sys_lin_2}).

\begin{prop}\label{controllability-prop}
For system (\ref{sys_lin_2}) to be $\kk(\delta)\left[\ddt\right]$-torsion free controllable (see \cite{Mounier_95, Rudolph_03}), it is necessary that $B$ and $F$, defined by (\ref{Fdef}), are $\kk(\delta)\left[\ddt\right]$-hyper-regular.

Moreover, in this case, there exists a polynomial $\bar{\pi}$ such that the localized system module at the powers of $\bar{\pi}$ is free.
\end{prop}
\begin{proof}
Assume that the system is $\kk(\delta)\left[\ddt\right]$ torsion free controllable and that $B$ is not $\kk(\delta)\left[\ddt\right]$-hyper-regular. Then, using the decomposition (\ref{Bdecomp}), $\Delta_B$ has at least one diagonal element which is a polynomial of degree larger than or equal to 1 with respect to $\ddt$ and with coefficients in $\kk(\delta)$. There indeed exists a non zero element $v$ of the $\kk[\delta,\ddt]$-free module generated by the components of $u$,  such that $\Delta_{B}v=0$ ($v$ is a non trivial solution of a differential delay equation). It is immediately seen that the pair $x=0, u=Nv$ is a non zero torsion element of the system module, which contradicts the freeness assumption.

 If $F$ is not $\kk(\delta)\left[\ddt\right]$-hyper-regular, its decomposition is given by $UF\tilde{Q} = \left( \Delta_{F}, 0 \right)$, $ \Delta_{F}$ having at least one diagonal element which is a polynomial of degree larger than or equal to 1 with respect to $\ddt$ and with coefficients in $\kk(\delta)$ which shows, using the representation (\ref{semi_impl_rep}),  that every  pair $(\xi_1,\xi_2)$ such that $\Delta_{F} \xi_1 =0$,  $\xi_1 \neq 0$, and $\xi_2$ arbitrary, satisfies  $\left( \Delta_{F}, 0 \right) \left( \begin{array}{c}\xi_1\\\xi_2\end{array}\right)=0$. Thus, the pair $(x,u)$ with $x=\tilde{Q}\left( \begin{array}{c}\xi_1\\\xi_2\end{array}\right)$ and $u$ satisfying $\Delta_{B}N^{-1}u=(I_{m},0_{m,n-m})MA\tilde{Q} \left( \begin{array}{c}\xi_1\\\xi_2\end{array}\right)$ (see (\ref{semi_impl_rep})) is a torsion element  of the system module. Consequently, the latter module cannot be $\kk(\delta)\left[\ddt\right]$-free.

To prove the existence of $\bar{\pi}$, we remark that, according to the Smith-Jacobson decomposition algorithm (see Section \ref{SJalgo-sec} of the Appendix), if $M\in \lsm{B}$ and $N\in \rsm{B}$ are such that $MBN=\left(\begin{array}{c}I_{m}\\0\end{array}\right)$, each row of $M$ and $N$ may contain the inverse of a polynomial of $\kk[\delta]$. Taking the LCM, say $\pi_{M,N}$, of these polynomials for all rows, we immediately get that $\pi_{M,N}\cdot M\in \uu{n}, \pi_{M,N}\cdot N\in \uu{m}$ with $\pi_{M,N} \in \kk[\delta]$. The same argument applies for $F$: consider $U\in\uud{n-m},\wt{Q}\in\uud{n}$ as above, namely such that $UF\wt{Q}=\left(I_{n-m}\ 0_{n-m,m}\right)$. There exists $\pi_{U,\tilde{Q}}\in \kk[\delta]$ such that $\pi_{U,\tilde{Q}}\cdot U \in \uu{n-m}$, $\pi_{U,\tilde{Q}}\cdot \tilde{Q} \in \uu{n}$. Taking $\bar{\pi}$ as the LCM of $\pi_{M,N}$ and $\pi_{U,\tilde{Q}}$, it is immediately seen from the decompositions of $B$ and $F$, multiplied by suitable powers of $\bar{\pi}$, that the system module over the localized ring $\kk[\delta,\bar{\pi}^{-1},\ddt]$ is free, and the proof is complete.
\end{proof}

\subsection{Differential $\pi$-Flatness}\label{pi-flat-sec}
We first recall the classical definition of a flat system \cite{Levine_09,Ramirez_04}, in the context of systems described by ordinary nonlinear differential equations: a system is said to be differentially flat if and only if there exists a set of independent variables, referred to as a flat output, such that every system variable (including the input variables) is a function of the flat output and a finite number of its successive time derivatives. More precisely, the system
\begin{equation}
\dot{x}=f(x,u)\nonumber
\end{equation}
with $x\in\RR^{n}$ and $u\in\RR^{m}$ is differentially flat if and only if there exist a set of independent variables (flat output)
\begin{equation}\label{y-flat}
y=h(x,u,\dot{u},\ddot{u},\ldots,u^{(r)}),\quad y\in\RR^{m}
\end{equation}
such that
\begin{eqnarray}
x&=&\alpha(y,\dot{y},\ddot{y},\ldots,y^{(s)})\\
u&=&\beta(y,\dot{y},\ddot{y},\ldots,y^{(s+1)})
\end{eqnarray}
and such that the system equations
\begin{equation}
\frac{d\alpha}{dt}(y,\dot{y},\ddot{y},\ldots,y^{(s+1)})=f\left(\alpha(y,\dot{y},\ddot{y},\ldots,y^{(s)}),\beta(y,\dot{y},\ddot{y},\ldots,y^{(s+1)})\right)
\end{equation}
are identically satisfied for all smooth enough function $t\mapsto y(t)$, $r$ and $s$ being suitable finite $m$-tuples of integers.

Let us now discuss further extensions of this definition to linear delay systems and consider first an elementary  example to motivate the next definition.

\begin{ex}
Let us consider the elementary system $$\dot{x}(t)=u(t-\tau)$$ with $A=\ddt$ and $B=\delta$ in the notations of (\ref{sys_lin_2}). Clearly, if we set $y=x$, $y$ looks like a flat output though $u=\delta^{-1}y$ contains a one-step prediction and belongs to $\kk(\delta)[\ddt]$ but not to $\kk[\delta,\ddt]$. However, for motion planning, such a dependence remains acceptable (see Remark~\ref{finitesupp-rem}), even if for feedback design it poses more delicate problems. This notion is called differential $\delta$-flatness (see e.g. \cite{Mounier_95,Rudolph_03}).
\end{ex}

According to \cite{Mounier_95}, this notion is generalized as follows:
\begin{defn}[Differential $\pi$-flatness \cite{Mounier_95}]\label{pi_flatness} The linear delay system (\ref{sys_lin_2}) is said to be differentially $\pi$-flat (or $\pi$-free) if and only if there exists a polynomial $\pi\in \kk[\delta]$, and a collection $y$ of $m$ $(\delta,\pi^{-1})$-differentially independent variables\footnote{more precisely, there does not exist a non zero matrix $S \in \mmp{m}{m}$ such that $Sy=0$,  or equivalently, $Sy=0$ implies $S\equiv 0$}, called $\pi$-flat output, of the form
\begin{equation}
y =P_{0}(\delta,\pi^{-1})x+P_{1}(\delta,\pi^{-1},\ddt)u,\label{y_pi_flat}
\end{equation}
with $P_{0}(\delta,\pi^{-1})\in(\kk[\delta,\pi^{-1}])^{m\times n}$ the set of matrices of size $m\times n$, with coefficients in $\kk[\delta,\pi^{-1}]$, $P_{1}(\delta,\pi^{-1},\ddt)\in\mmp{m}{m}$, and such that
\begin{eqnarray}
x &=&Q(\delta,\pi^{-1},\ddt)y,\label{x_pi_flat}\\
u &=&R(\delta,\pi^{-1},\ddt)y,\label{u_pi_flat}
\end{eqnarray}
with $Q(\delta,\pi^{-1},\ddt)\in\mmp{n}{m}$, and $R(\delta,\pi^{-1},\ddt)\in\mmp{m}{m}$.
\end{defn}
In other words, definition \ref{pi_flatness} states that the components of a $\pi$-flat output $y$ can be obtained as a $\kk[\delta,\pi^{-1},\ddt]$-linear combination of the system variables, and that the system variables $(x, u)$ are also $\kk[\delta,\pi^{-1},\ddt]$-linear combinations of the components of $y$. Thus $x$ and $u$ can be calculated from $y$ using differentiations, delays, and predictions (coming from the inverse of $\pi$). Note that the facts that every element of the system module is a $\kk[\delta,\pi^{-1},\ddt]$-linear combination of the components of $y$, and that the components of $y$ are $\kk[\delta,\pi^{-1},\ddt]$-independent, indeed imply that the components of $y$ form a basis of the system module $\kk[\delta,\pi^{-1},\ddt] \otimes_{\kk[\delta,\ddt]} \Lambda$, which is therefore free, hence the equivalence with the definition of subsection~\ref{module-sec}.

If system (\ref{sys_lin_2}) is considered in implicit form (\ref{impl_rep}) after elimination of the input $u$, since this elimination expresses $u$ as a $\kk\left(\delta\right)[\ddt]$-combination of the components of $x$, the expression (\ref{y_pi_flat}), combined with (\ref{udef}), reads $y=Px$ with $P\in\mmd{m}{n}$.
The previous definition is thus adapted as follows:
\begin{defn}\label{pi_flatness_impl} The implicit linear delay system (\ref{impl_rep}) is said to be differentially $\pi$-flat if and only if there exists a polynomial $\pi\in \kk[\delta]$, and a collection $y$ of $m$ $(\delta,\pi^{-1})$-differentially independent variables, called $\pi$-flat output, of the form
\begin{equation}
y =P(\delta,\pi^{-1},\ddt)x,\label{y_pi_flat_impl}
\end{equation}
with $P(\delta,\pi^{-1},\ddt)\in\mmp{m}{n}$, and such that
\begin{equation}
x =Q(\delta,\pi^{-1},\ddt)y,\label{x_pi_flat_impl}
\end{equation}
with $Q(\delta,\pi^{-1},\ddt)\in\mmp{n}{m}$.
\end{defn}

The matrices $P$, $Q$ and $R$ of (\ref{y_pi_flat})--(\ref{u_pi_flat}) in the explicit case, and $P$ and $Q$ of (\ref{y_pi_flat_impl})--(\ref{x_pi_flat_impl}) in the implicit case, are called \emph{defining operators} of the $\pi$-flat output $y$.

\begin{rem}
In \cite{Mounier_95} and later (see e.g. \cite{Rudolph_03}), the above notion is often called $\pi$-freeness and introduced via the notion of system module. The wording $\pi$-flatness appears, to the authors knowledge, for the first time in \cite{Petit_00}.  It has also been related to system parameterization in \cite{Pommaret_99,Chyzak_07}. We have preferred here the name $\pi$-flatness, in reference to differential flatness, and to directly present it via the notion of flat output, rather than basis  of the system module. Note that in formula~(\ref{y_pi_flat}), $P_{0}$ is a 0th degree polynomial of $\ddt$ to mimic the general definition in (\ref{y-flat}) that does not include time derivatives of $x$, with (\ref{y_pi_flat})--(\ref{u_pi_flat}) restricted to linear expressions.
\end{rem}

\begin{ex}\label{ex_contd}
Let us go back again to example 1 and let us prove that $y=x_{1}$ is a $\pi$-flat output with $\pi=(1-\delta)\delta^{2}$. From (\ref{example_1_bis}), we have 
$$x_{2}=\delta^{-1}(1-\delta)^{-1}\frac{1}{k}\dot{x}_{1}=\delta^{-1}(1-\delta)^{-1}\frac{1}{k}\dot{y}$$
\begin{equation}
u=\delta^{-1}\ddt x_{2}=\delta^{-2}(1-\delta)^{-1}\left(-\frac{\dot{k}}{k^{2}}\dot{y}+\frac{1}{k}\ddot{y}\right)\label{u_ex1}
\end{equation}
In other words, following the notations of (\ref{y_pi_flat})--(\ref{u_pi_flat}), $P_{0}=\left( 1,0\right)$, $P_{1}=0$ and
\begin{equation}
\left(\begin{array}{c}x_{1}\\x_{2}\end{array}\right)=\left(\begin{array}{c}1\\\pi^{-1}\delta\frac{1}{k}\ddt\end{array}\right)y \triangleq Q(\delta,\pi^{-1},\ddt)y\label{x_ex1}
\end{equation}
\begin{equation}
u=\pi^{-1}\left(- \frac{\dot{k}}{k^{2}}\ddt+\frac{1}{k}\frac{d^{2}}{dt^{2}}\right)y \triangleq R(\delta,\pi^{-1},\ddt)y\label{u_bis_ex1}
\end{equation}
which proves that $y$ is a $\pi$-flat output.
\end{ex}

\section{Main Result}\label{main-sec}

In this section, we propose a simple and effective algorithm for the computation of $\pi$-flat outputs of linear time-delay systems based on the following necessary and sufficient condition for the existence of defining operators of a $\pi$-flat output. Moreover, explicit expressions of $\pi$ and of these operators are obtained.

\begin{thm}\label{pi_flat_thm}
A necessary and sufficient condition for system (\ref{sys_lin_2}) to be $\pi$-flat is that the matrices $B$ and $F$ are $\kk(\delta)\left[\ddt\right]$-hyper-regular. 

We construct the operators $P$, $Q$ and $R$  and the polynomial $\pi$  as follows.
\begin{itemize}
\item[\bf{0.}]  According to Propositions~\ref{impl-prop} and \ref{controllability-prop}, construct the decomposition of $B$ (\ref{Bdecomp}), define $F$ by (\ref{Fdef}) and compute $\bar{\pi}$;
\item[\bf{1.}]  $Q=\wt{Q}\left( \begin{array}{c}0_{n-m,m}\\I_{m}\end{array}\right),~\text{with}~\wt{Q}\in\rsm{F}$. Note that  $\bar{\pi}\cdot Q \in \mm{n}{m}$ by construction;
\item[\bf{2.}]  $R=N(I_{m},0_{m,n-m})MAQ,~\text{with}~N\in\rsm{B}$. There exists $\pi_{R} \in \kk[\delta]$ such that $\pi_{R}\cdot R\in \mm{m}{m}$;
\item[\bf{3.}]  $P=W(I_{m},0_{m,n-m})\wt{P},~\text{with}~\wt{P}\in\lsm{Q}~\text{and}~W\in\rsm{Q}$. There exists $\pi_{P}\in \kk[\delta]$ such that $\pi_{P}\cdot P\in \mm{m}{n}$.
\end{itemize}
Let $\pi$ be given by $\pi= \lcm{\bar{\pi},\pi_{P},\pi_{R}}$, the least common multiple of $\bar{\pi}$, $\pi_{P}$ and $\pi_{R}$.

Thus, $P$, $Q$ and $R$ are defining operators with $y=Px$, $x=Qy$, $u=Ry$, and $y$ is a $\pi$-flat output. 
\end{thm}

\begin{proof}
If $B$ or $F$ are not $\kk(\delta)\left[\ddt\right]$-hyper-regular, according to Proposition~\ref{controllability-prop},  the $\kk(\delta)[\ddt]$ system module cannot be torsion free. Therefore, system (\ref{sys_lin_2}) cannot be  $\pi$-flat, for any $\pi\in \kk[\delta]$. Taking the contrapositive of this statement, the hyper-regularity of $B$ and $F$ is proven to be necessary.

We now prove that the $\kk(\delta)\left[\ddt\right]$-hyper-regularity of $B$ and $F$ is sufficient to construct the defining matrices $P$, $Q$ and $R$ as well as a liberation polynomial $\pi$.

We first use the implicit form (\ref{impl_rep}) of Proposition~\ref{impl-prop}, and more precisely (\ref{Fdef})--(\ref{udef}) to obtain a decomposition of $B$ and $F$.
Since $F$ is hyper-regular by assumption, there exist $U\in\uud{n-m},\wt{Q}\in\uud{n}$ such that $UF\wt{Q}=\left(I_{n-m}\ 0_{n-m,m}\right)$. Consequently $Q=\wt{Q}\left(\begin{array}{c}0_{n-m,m}\\\quad I_{m}\end{array}\right)$, with $\wt{Q}\in \rsm{F}$, is such that $FQ=0_{n-m,m}$. Thus setting $x=Qy$ we have that $FQy=0$ for all $y$.

The existence of $\bar{\pi}\in \kk[\delta]$ is proved following the same argument as in the proof of Proposition~\ref{controllability-prop}: for an arbitrary $\alpha \times \beta$ hyper-regular matrix $M$, according to the Smith-Jacobson decomposition algorithm (see Section \ref{SJalgo-sec} of the Appendix), if $U_{M}\in \lsm{M}$ and $V_{M}\in \rsm{M}$ are such that $U_{M}MV_{M}=(I_{\alpha},0)$ if $\alpha\leq \beta$ (resp.
$U_{M}MV_{M}=\left(\begin{array}{c}I_{\beta}\\0\end{array}\right)$ if $\alpha\geq \beta$), each row of $U_{M}$ and $V_{M}$ may contain the inverse of a polynomial of $\kk[\delta]$. Taking the LCM, say $\pi_{M}$, of these polynomials for all rows, we immediately get that $\pi_{M}\cdot U_{M}\in \mm{\alpha}{\beta}$ with $\pi_{M} \in \kk[\delta]$. Applying this result to the decompositions of $B$ and $F$, we have proven the existence of $\bar{\pi}$ such that $\bar{\pi}\cdot \wt{Q} \in \uu{n}$ and $\bar{\pi}\cdot Q\in \mm{n}{m}$, which proves item 1.

Going back to (\ref{Bdecomp}) and (\ref{udef}), setting $R=N\left(I_{m}\quad 0_{m,n-m}\right)MAQ$, we obtain $u=Ry$ with $N\in\rsm{B}$. Finally, the proof of the existence of $\pi_{R} \in \kk[\delta]$ such that $\pi_{R}\cdot R\in \mm{m}{m}$ follows the same lines as in item~1, which proves 2.

Since $Q$ is hyper-regular by construction, its Smith-Jacobson decomposition yields the existence of $\wt{P}\in \uud{n}$ and $W\in \uud{m}$ such that 
$$\wt{P}QW=\left(\begin{array}{c}I_{m}\\0_{n-m,m}\end{array}\right)$$ 
thus $\left(I_{m}\quad 0_{m,n-m}\right)\wt{P}QW=I_{m}$ and, setting 
$$P=W(I_{m},0_{m,n-m})\wt{P}$$ 
it results that $Px=W(I_{m},0_{m,n-m})\wt{P}x=W(I_{m},0_{m,n-m})\wt{P}Qy=WW^{-1}y=y$, and the proof of the third item is complete, noting again that the existence of $\pi_{P} \in \kk[\delta]$ such that $\pi_{P}\cdot P\in \mm{m}{n}$ follows the same lines as in  items 1 and 2. 

Finally, taking $\pi= \lcm{\bar{\pi},\pi_{P},\pi_{R}}$, the least common multiple of $\bar{\pi}$, $\pi_{P}$ and $\pi_{R}$, it is straightforward to show that $\pi \cdot P \in \mm{m}{n}$, $\pi\cdot  Q \in\mm{n}{m}$ and $\pi \cdot R \in\mm{m}{m}$. Therefore, $P$, $Q$ and $R$ are defining matrices of a $\pi$-flat output for system (\ref{sys_lin_2}) and thus that system (\ref{sys_lin_2}) is $\pi$-flat.   
\end{proof}

Theorem \ref{pi_flat_thm} is easily translated into the Algorithm~\ref{algo_pi_flat} presented below.

\begin{algorithm}[h]\label{algo_pi_flat}
\SetAlgoNoLine
\LinesNumbered
\DontPrintSemicolon
\SetKwInput{KwInit}{Initialization}
\SetKwInput{Algo}{Algorithm}
\caption{Procedure to compute $\pi$-flat outputs.}
\KwIn{Two matrices $A\in\mm{n}{n}$ and $B\in\mm{n}{m}$.}
\KwOut{a polynomial $\pi\in \kk[\delta]$ and defining operators $P\in\mmp{m}{n}$, $Q\in\mmp{n}{m}$ and $R\in\mmp{m}{m}$ such that $y=Px$, $x=Qy$ and $u=Ry$.}
\KwInit{Test of hyper-regularity of $B$ by its Smith-Jacobson decomposition, which also provides $M\in\lsm{B}$ and $N\in\rsm{B}$, i.e. such that $MBN=\left( I_{m}\, ,\, 0_{n-m,m}\right)^{T}$. If $B$ is not hyper-regular, the system is not $\pi$-flat whatever $\pi\in \kk[\delta]$.} 
\Algo

Set $F=(0_{n-m,m},I_{n-m})MA$ and test if $F$ is hyper-regular by computing its Smith-Jacobson decomposition: $VF\wt{Q}=(I_{n-m},0_{n-m,m})$. If $F$ is not hyper-regular, the system is not $\pi$-flat whatever $\pi\in \kk[\delta]$.  Otherwise, compute $\bar{\pi}\in\kk[\delta]$ such that $\bar{\pi}\cdot M\in \uu{n}$, $\bar{\pi}\cdot N\in \uu{m}$, $\bar{\pi}\cdot \wt{Q}\in \uu{n}$.\;
Compute $Q=\wt{Q}\left( 0_{n-m,m}\, ,\, I_{m} \right)^{T}$,  $R=N(I_{m},0_{m,n-m})MAQ$ and $\pi_{R}\in \kk[\delta]$ such that $\pi_{R}\cdot R\in \mm{m}{m}$.\;
Compute a Smith-Jacobson decomposition of $Q$: $\wt{P}QW=\left( I_{m}\, ,\, 0_{n-m,m}\right)^{T}$, $P=W(I_{m},0_{m,n-m})\wt{P}$, and a polynomial $\pi_{P} \in\kk[\delta]$ such that $\pi_{P}\cdot P \in \mm{m}{n}$.\;
Compute the polynomial $\pi=\lcm{\bar{\pi}, \pi_{P}, \pi_{R}} \in\kk[\delta]$. The system is $\pi$-flat.
\end{algorithm}

\begin{rem}
The $\pi$-flatness criterion is given in terms of properties of the matrices $B$ and $F$ that depend only on the larger ring $\kk(\delta)[\ddt]$. If, in addition, the submodule of $\kk[\delta]$ generated by the powers of  $\pi$ is torsion free (e.g. $\pi=\delta$, for which the equation $\delta f=0$ admits the unique solution $f=0$), then the system module $\Lambda$ over $\kk[\delta,\pi^{-1},\ddt]$ is torsion free, and the original module $\Lambda$ (over $\kk[\delta,\ddt]$) is free if, and only if, $\pi =1$. Note that computing a free basis of the system module directly over $\kk[\delta,\ddt]$, would require more elaborated tools such as those developed in \cite{Pommaret_99,Pommaret_01,Quadrat_07}.
\end{rem}

\begin{rem}\label{multdelay-rem}
Since the computations are made in the larger modules $\mmd{p}{q}$ for suitable $p$ and $q$, Theorem~\ref{pi_flat_thm} remains valid for systems depending on an arbitrary but finite number of delays, say $\delta_{1},\ldots, \delta_{s}$, by replacing the field $\kk(\delta)$, on which $\mmd{p}{q}$ is modeled, by the fraction field $\kk(\delta_{1},\ldots, \delta_{s})$ generated by the multivariate polynomials of $\delta_{1},\ldots, \delta_{s}$. An example with two independent delays is presented in Example~\ref{string-ex-sec}.
\end{rem}

\begin{rem}
Our results also apply in the particular case of linear time-varying systems without delays. If $\kk$ is the field of meromorphic functions of time, $\kk[\ddt]$ is a Principal Ideal Domain and all the computations involved in Theorem~\ref{pi_flat_thm} and the associated algorithm, remain in this ring, contrarily to the case with delay. Therefore, the last step consisting in finding the so-called liberation polynomial $\pi$ is needless. Related results may be found in \cite{Levine_03,Trentelman_04,Pommaret_99}.
\end{rem}
\section{Examples}\label{ex-sec}

\subsection{Back to the Introductory Example}
Going back to Example~\ref{Example_1}, let us apply the previous algorithm to the time-delay system defined by (\ref{example_1}), for which a $\pi$-flat output is already known from Example \ref{ex_contd}. The first step, consisting in the computation of a Smith-Jacobson decomposition of the matrix $B$ has already been done, the left and right unimodular matrices $M\in \lsm{B}$ and $N\in \rsm{B}$ being given by (\ref{Smith_B}) with $M=U'$ and $N=V=1$. Then the matrix $F=(0\quad 1)MA$ of an implicit representation of (\ref{example_1}) is given by (\ref{implicit_form}). Its Smith-Jacobson decomposition $VF\wt{Q}=(1\quad 0)$ is given by (\ref{F_decomp})-(\ref{Smith_VF}), with $V=U_{F}=1$ and $\wt{Q}=V_{F}$, and has been seen to be hyper-regular in Example~\ref{example_smith_bis}. We easily check that $\bar{\pi}= (1-\delta)\delta$.

Going on with step 2, we set
$$Q=\wt{Q}\left(\begin{array}{c}0\\1\end{array}\right)=\left(\begin{array}{c}1\\\ds \delta^{-1}(1-\delta)^{-1} \frac{1}{k(t)}\ddt\end{array}\right)$$ 
and verify that  $\bar{\pi}\cdot Q \in \mm{2}{1}$. Moreover
$$R=N(1\quad 0)MAQ=\delta^{-2}(1-\delta)^{-1}\left(-\frac{\dot{k}(t)}{k^{2}(t)}\ddt+\frac{1}{k(t)}\frac{d^{2}}{dt^{2}}\right)$$
and we have $\pi_{R}=\bar{\pi}\delta$.

Note that, setting $x=Qy$ and $u=Ry$, we recover formulae (\ref{x_ex1}) and (\ref{u_ex1}) (or equivalently (\ref{u_bis_ex1})). 

According to step 3 of the algorithm, we compute a Smith-Jacobson decomposition of $Q$: $\wt{P}QW=\left(\begin{array}{c}1\\ 0\end{array}\right)$, which provides $W=1$ and
$$\wt{P}=\left(\begin{array}{cc}1&0\\-\delta^{-1}(1-\delta)^{-1}\frac{1}{k(t)}\ddt&1\end{array}\right)$$
hence $P=\left(1\quad 0\right)$ and $y=Px=x_{1}$. Here, $\pi_{P}=1$.
Finally, the least common multiple of $(1-\delta)\delta$, $(1-\delta)\delta^{2}$ and 1 is $\pi=(1-\delta)\delta^{2}$. We have thus verified that Algorithm~\ref{algo_pi_flat} comes up with the same conclusion as Example~\ref{ex_contd}.

\subsection{A Multi-input Example}
Let us consider the following academic example of multi-input delay system:
\small
\begin{equation}
\left\{\begin{array}{lll}
\dot{x}_{1}(t)+x_{1}^{(2)}(t)-2x_{1}^{(2)}(t-\tau)+x_{1}^{(3)}(t)+x_{1}^{(4)}(t-\tau)-x_{2}^{(3)}(t)+x_{2}^{(5)}(t)-x_{3}^{(2)}(t)\\
-\dot{x}_{4}(t)+x_{4}^{(3)}(t)=u_{1}(t)+\dot{u}_{1}(t)+u_{2}(t),\\
\dot{x}_{1}(t)+\dot{x}_{1}(t-\tau)-\dot{x}_{1}(t-2\tau)+x_{1}^{(2)}(t)+x_{1}^{(2)}(t-\tau)+x_{1}^{(2)}(t-2\tau)-x_{1}^{(3)}(t-\tau)\\
+2\dot{x}_{2}(t)+\dot{x}_{2}(t-\tau)-x_{2}^{(2)}(t)-x_{2}^{(4)}(t)+\dot{x}_{3}(t)+x_{3}^{(2)}(t-\tau)-x_{4}(t)-x_{4}(t-\tau)\\
-x_{4}^{(2)}(t)=u_{1}^{(2)}(t-\tau)+\dot{u}_{2}(t-\tau),\\
-x_{1}(t-2\tau)+\dot{x}_{1}(t-3\tau)+x_{1}^{(2)}(t-2\tau)-x_{2}(t-\tau)+\dot{x}_{2}(t-2\tau)+x_{2}^{(3)}(t-\tau)\\
-x_{3}(t-\tau)+\dot{x}_{3}(t-2\tau)+\dot{x}_{4}(t-\tau)=\dot{u}_{1}(t-2\tau)+u_{2}(t-2\tau),\\
\dot{x}_{1}(t-\tau)+\dot{x}_{2}(t)+\dot{x}_{3}(t)=\dot{u}_{1}(t)+u_{2}(t).\\
\end{array}\right.\label{sys_ex_2}
\end{equation}
\normalsize
Denoting by $x$ the state vector, $x=(x_{1},x_{2},x_{3},x_{4})^{T}$, by $u$ the input vector, $u=(u_{1},u_{2})^{T}$, and by $\delta$ the delay operator of length $\tau$, system (\ref{sys_ex_2}) can be rewritten in matrix form $Ax=Bu$, with $A\in\mm{4}{4}$ and $B\in\mm{4}{2}$ defined by
\begin{equation}
A=\left(\begin{array}{cccc}A_{1}&A_{2}&A_{3}&A_{4}\end{array}\right)
\end{equation}
with 
\begin{equation}
\begin{array}{l}
A_{1}=\left( \begin{array}{c}
\frac{d}{dt}+\frac{{d}^{2}}{{dt}^{2}}\left( 1-2\delta \right) +\frac{{d}^{3}}{{dt}^{3}}+\frac{{d}^{4}}{{dt}^{4}}\delta\\
{\frac{d}{dt}}\left( 1+ \delta - \delta^{2} \right) + {\frac{{d}^{2}}{{dt}^{2}}}\left( 1+ \delta + \delta^{2} \right) - {\frac{{d}^{3}}{{dt}^{3}}}\delta
\\
- {\delta}^{2} + {\frac{d}{dt}}{\delta}^{3} + {\frac{{d}^{2}}{{dt}^{2}}{\delta}^{2}}
\\
{\frac{d}{dt}\delta}
\end{array}\right)
\\
A_{2}=\left( \begin{array}{c}
-{\frac{{d}^{3}}{{dt}^{3}}}+{\frac{{d}^{5}}{{dt}^{5}}}
\\
\,{\frac{d}{dt}}\left( 2 + \delta \right) - {\frac{{d}^{2}}{{dt}^{2}}} - {\frac{{d}^{4}}{{dt}^{4}}}
\\
- \delta + {\frac{d}{dt}{\delta}^{2}} + {\frac{{d}^{3}}{{dt}^{3}}\delta}
\\
{\frac{d}{dt}}
\end{array}\right)
\\
A_{3}=\left( \begin{array}{c}
-{\frac{{d}^{2}}{{dt}^{2}}}
\\
{\frac{d}{dt}}+{\frac{{d}^{2}}{{dt}^{2}}\delta}
\\
- \delta + {\frac{d}{dt}{\delta}^{2}}
\\
{\frac{d}{dt}}
\end{array}\right),
\quad
A_{4}=\left( \begin{array}{c}
-{\frac{d}{dt}}+{\frac{{d}^{3}}{{dt}^{3}}}
\\
-\left(1+ \delta \right) - {\frac{{d}^{2}}{{dt}^{2}}}
\\
 {\frac{d}{dt}\delta}
 \\
0
\end{array}\right)
\end{array}
\end{equation}
and
\begin{equation}
B=\left(\begin{array}{cc}
1+ \frac{d}{dt}
&
1
\\
\frac{{d}^{2}}{{dt}^{2}}\delta
&
\frac{d}{dt}\delta
\\
\frac{d}{dt}{\delta}^{2}
&
{\delta}^{2}
\\
\frac{d}{dt}
&
1
\end{array}\right).\nonumber
\end{equation}
We apply Algorithm~\ref{algo_pi_flat} to compute a $\pi$-flat output if it exists. We start with the Smith-Jacobson decomposition of $B$. By left multiplying $B$ by the following product of unimodular matrices $M=M_{4}M_{3}M_{2}M_{1} \in\uud{4}$, given by
\begin{equation}
\begin{array}{ll}
M_{1}=\left(\begin{array}{cccc}1&0&0&-1\\0&1&0&0\\0&0&1&0\\0&0&0&1\end{array}\right),&
M_{2}=\left(\begin{array}{cccc}1&0&0&0\\-\frac{d^{2}}{dt^{2}}\delta&1&0&0\\-\frac{d}{dt}\,{\delta}^{2}&0&1&0\\-\frac{d}{dt}&0&0&1\end{array}\right),\\\\
M_{3}=\left(\begin{array}{cccc}1&0&0&0\\0&0&0&1\\0&1&0&0\\0&0&1&0\end{array}\right),&
M_{4}=\left(\begin{array}{cccc}1&0&0&0\\0&1&0&0\\0&-\frac{d}{dt}\delta&1&0\\0&-{\delta}^{2}&0&1\end{array}\right),
\end{array}\nonumber
\end{equation}
and by setting $N=I_{2}$, we obtain the Smith-Jacobson decomposition
\begin{equation}
MBN=\left(
\begin{array}{cccc}
1&0&0&-1\\-\frac{d}{dt}&0&0&\frac{d}{dt}+1\\0&1&0&-\frac{d}{dt}\delta\\0&0&1&-\delta^{2}\end{array}\right)
\left(\begin{array}{cc} 
\frac{d}{dt}+1&1\\\frac{{d}^{2}}{{dt}^{2}}\delta&\frac{d}{dt}\delta\\\frac{d}{dt}{\delta}^{2}&{\delta}^{2}\\\frac{d}{dt}&1
\end{array}\right)=\left(
\begin{array}{cccc}1&0\\0&1\\0&0\\0&0\end{array}\right),\nonumber
\end{equation}
thus showing that $B$ is hyper-regular. We then compute an implicit representation of (\ref{sys_ex_2}) by
$F=(0_{2,2}\quad I_{2})MA$, i.e.
\begin{equation}
F=\left(\begin{array}{cccc} F_{1}&F_{2}&F_{3}&F_{4}\end{array}\right)
\end{equation}
with
\begin{equation}
\begin{array}{l}
F_{1}=\left(
\begin{array}{c}
{\frac{d}{dt}}\left(1+ \delta-\delta^{2}\right) + {\frac{{d}^{2}}{{dt}^{2}}}\left( 1+ \delta \right) -{\frac{{d}^{3}}{{dt}^{3}}\delta}
\\
-{\delta}^{2} + {\frac{{d}^{2}}{{dt}^{2}}{\delta}^{2}}
\end{array}
\right)
\\
F_{2}=\left(
\begin{array}{c}
{\frac{d}{dt}}\left( 2 + \delta \right) - {\frac{{d}^{2}}{{{dt}}^{2}}} \left( 1 + \delta \right) - {\frac{{d}^{4}}{{dt}^{4}}}
\\
-\delta+{\frac{{d}^{3}}{{dt}^{3}}\delta}
\end{array}
\right)
\\
F_{3}=\left(
\begin{array}{c}
{\frac{d}{dt}}
\\
-\delta
\end{array}
\right),
\quad
F_{4}=\left(
\begin{array}{c}
-1-\delta-{\frac{{d}^{2}}{{dt}^{2}}}
\\
{\frac{d}{dt}\delta}
\end{array}
\right)
\end{array}\
\end{equation}
to which corresponds the difference-differential system
\begin{equation}
\begin{array}{l}
\ds \left( \dot{x}_{1}(t)+\dot{x}_{1}(t-\tau)-\dot{x}_{1}(t-2\tau)+\ddot{x}_{1}(t)+\ddot{x}_{1}(t-\tau)-x_{1}^{(3)}(t-\tau)\right) \\
\ds \hspace{2cm} +\left( 2\dot{x}_{2}(t)+\dot{x}_{2}(t-\tau) -\ddot{x}_{2}(t)-\ddot{x}_{2}(t-\tau)-x_{2}^{(4)}(t)\right)  \vspace{0.4em}\\
\ds \hspace{3cm} + \dot{x}_{3}(t)- \left( x_{4}(t) + x_{4}(t-\tau) + \ddot{x}_{4}(t)\right) =0,\vspace{0.4em}\\
\ds \left( -x_{1}(t-2\tau)+\ddot{x}_{1}(t-2\tau)\right) - \left( x_{2}(t-\tau)- x_{2}^{(3)}(t-\tau)\right) \vspace{0.4em}\\
\ds \hspace{3cm} - x_{3}(t-\tau) + \dot{x}_{4}(t-\tau)=0.
\end{array}\nonumber
\end{equation}
According to step 1, we compute a right Smith-Jacobson decomposition of $F$:
\begin{equation}
VF\wt{Q}=\left(\begin{array}{cccc}1&0&0&0\\0&1&0&0\end{array}\right)
\end{equation}
where $V=1$ and 
\begin{equation}\label{Qtilde-ex2}
\wt{Q}=\left(\begin{array}{cccc} 0&0&0&1\\0&0&1&0\\-\left(1+\delta\right)^{-1}{\frac{d}{dt}}&-{\delta}^{-1}\left(1+\delta\right)^{-1}\left(1+\delta+{\frac{{d}^{2}}{{dt}^{2}}}\right)&{\frac{{d}^{2}}{{dt}^{2}}}-1&{\frac{{d}^{3}}{{dt}^{3}}}+{\frac{{d}^{2}}{{dt}^{2}}}-\delta\\-\left(1+\delta\right)^{-1}&-\delta^{-1}\left(1+\delta\right)^{-1}{\frac{d}{{dt}}}&\left(1-{\frac{d}{{dt}}}\right)\frac{d}{dt}&\left({\frac{d}{dt}}+1-\delta\right)\frac{d}{dt}\end{array}\right),\nonumber
\end{equation} 
showing thus that $F$ is hyper-regular and that $\bar{\pi}=\delta(1+\delta)$.

For the interested reader, $\wt{Q}$ is obtained as the product 
$\wt{Q}=\wt{Q}_{1}\wt{Q}_{2}\wt{Q}_{3}$ of matrices of elementary actions:
\begin{equation}
\wt{Q}_{1}=\left(\begin{array}{cccc}0&0&0&1\\0&0&1&0\\{\frac{d}{dt}}&1&0&0\\1&0&0&0\end{array}\right),\nonumber
\end{equation}
\begin{equation}
\wt{Q}_{2}=\left(\begin{array}{cccc}\wt{q}_{1,1,2}&\wt{q}_{1,2,2}&\wt{q}_{1,3,2}&\wt{q}_{1,4,2}
\\0&1&0&0\\0&0&1&0\\0&0&0&1\end{array}\right),\nonumber
\end{equation}
with
\begin{equation}
\begin{array}{l}
\wt{q}_{1,1,2} = -\left(1+\delta\right)^{-1}\\
\wt{q}_{1,2,2} = \left(1+\delta\right)^{-1}\frac{d}{dt}\\
\wt{q}_{1,3,2} = \left(1 + \delta \right)^{-1} \left( \left( 2+\delta \right) \frac{d}{dt} - \frac{d^{4}}{dt^{4}}\right) - \frac{d^{2}}{dt^{2}} \\
\wt{q}_{1,4,2} =\ddt - \delta^{2} \left(1 + \delta \right)^{-1}\ddt + \frac{d^{2}}{dt^{2}} - \left(1 + \delta \right)^{-1} \frac{d^{3}}{dt^{3}}\delta
\end{array}
\end{equation}
and
\begin{equation}
\wt{Q}_{3}=\left(\begin{array}{cccc}1&0&0&0\\0&-{\delta}^{-1}&-1+{\frac{{d}^{3}}{{dt}^{3}}}&{\frac{{d}^{2}}{{dt}^{2}}{\delta}}-{\delta}\\0&0&1&0\\0&0&0&1\end{array}\right).\nonumber
\end{equation}

According to step 2, we get
$$Q=\wt{Q}\left( \begin{array}{c} 0_{2,2}\\ I_{2}\end{array}\right)
=
\left( \begin{array}{cc}
0
&
1
\\
1
&
0
\\
{\frac{{d}^{2}}{{dt}^{2}}}-1
&
{\frac{{d}^{3}}{{dt}^{3}}}+{\frac{{d}^{2}}{{dt}^{2}}}-\delta
\\
\left(1-{\frac{d}{{dt}}}\right)\frac{d}{dt}&\left({\frac{d}{dt}}+1-\delta\right)\frac{d}{dt}
\end{array}\right)$$

\begin{equation}
R=N(I_{2}\quad 0_{2,2})MAQ=\left(\begin{array}{cc}-{\frac{{d}^{3}}{{dt}^{3}}}
&
{\frac{d}{dt}}-{\frac{{d}^{3}}{{{dt}}^{3}}}-{\frac{{d}^{4}}{{dt}^{4}}}
\\
\left({\frac{d}{dt}}+1\right)\frac{{d}^{3}}{{dt}^{3}}
&
-{\frac{{d}^{2}}{{dt}^{2}}}+{\frac{{d}^{3}}{{dt}^{3}}}+2\,{\frac{{d}^{4}}{{dt}^{4}}}+{\frac{{d}^{5}}{{dt}^{5}}}\end{array}\right),\nonumber
\end{equation}
with $\pi_{R}=1$.
From $x=Qy$ and $u=Ry$, we deduce the expressions
\begin{equation}\label{Q-ex2}
\left(\begin{array}{cc}x_{1}(t)\\x_{2}(t)\\x_{3}(t)\\x_{4}(t)\end{array}\right) = \left(\begin{array}{cc}
y_{2}(t)\\
y_{1}(t)\\-y_{1}(t)+y_{1}^{(2)}(t)-y_{2}(t-\tau)+y_{2}^{(2)}(t)+y_{2}^{(3)}(t)\\\dot{y}_{1}(t)-y_{1}^{(2)}(t)+\dot{y}_{2}(t)-\dot{y}_{2}(t-\tau)+y_{2}^{(2)}(t)
\end{array}\right).
\end{equation}
and
\begin{equation}\label{R-ex2}
\left(\begin{array}{cc}u_{1}(t)\\u_{2}(t)\end{array}\right) = \left(\begin{array}{cc}
-y_{1}^{(3)}(t)+\dot{y}_{2}(t)-y_{2}^{(3)}(t)-y_{2}^{(4)}(t)\\
y_{1}^{(3)}(t)+y_{1}^{(4)}(t)-y_{2}^{(2)}(t)+y_{2}^{(3)}(t)+2y_{2}^{(4)}(t)+y_{2}^{(5)}(t)
\end{array}\right).
\end{equation}

Next, according to step 3, we compute $\wt{P}\in\uud{4}$ and $W\in \uud{2}$ such that $\wt{P}QW=(I_{2}\quad 0_{2,2})^{T}$. 

\begin{equation}
\wt{P}=\left(\begin{array}{cccc}0&1&0&0\\1&0&0&0\\-{\frac{{d}^{3}}{{dt}^{3}}}-{\frac{{d}^{2}}{{dt}^{2}}}+\delta&-{\frac{{d}^{2}}{{dt}^{2}}}+1&1&0\\{\frac{d\delta}{dt}}-{\frac{{d}^{2}}{{dt}^{2}}}-{\frac{d}{dt}}&-{\frac{d}{dt}}+{\frac{{d}^{2}}{{dt}^{2}}}&0&1\end{array}\right),\quad W=I_{2}.\nonumber
\end{equation}
Again, $\wt{P}$ is obtained as the product $\wt{P}=\wt{P}_{2}\wt{P}_{1}$ of elementary actions:
\begin{equation}
\wt{P}_{2}=\left(\begin{array}{cccc}1&0&0&0\\0&1&0&0\\-{\frac{{d}^{2}}{{dt}^{2}}}+1&-{\frac{{d}^{3}}{{dt}^{3}}}-{\frac{{d}^{2}}{{dt}^{2}}}+\delta&1&0\\-\left(1-{\frac{d}{dt}}\right)\frac{d}{dt}&-\left({\frac{d}{dt}}+1-\delta\right)\frac{d}{dt}&0&1\end{array}\right) , \quad
\wt{P}_{1}=\left(\begin{array}{cccc}0&1&0&0\\1&0&0&0\\0&0&1&0\\0&0&0&1\end{array}\right).
\nonumber
\end{equation}

Then 
$$P=W\left(I_{2}\quad 0_{2,2}\right)\wt{P}= \left( \begin{array}{cccc}0&1&0&0\\1&0&0&0\end{array}\right).
$$
which, with $y=Px$, yields 
\begin{equation}\label{P-ex2}
y_{1}=x_{2}, \quad y_{2}=x_{1}
\end{equation}
and $\pi_{P}=1$.

Then it is immediately seen that $\pi=\delta(1+\delta)$ and that the system is $\pi$-flat.

\begin{rem}
It is worth noting that the polynomial $\pi=\delta(1+\delta)$ only appears at the intermediate level of the computation of $\wt{Q}$, and not anymore in the defining matrices $P$, $Q$ and $R$. However, this means that the system module contains elements $z$ such that $z(t)=-z(t-\tau)$ for all $t$, that satisfy $\pi z=0$, thus preventing this module from being free.
\end{rem}

\subsection{Vibrating String With an Interior Mass}\label{string-ex-sec}
As noted in Remark~\ref{multdelay-rem}, the computation of $\pi$-flat outputs based on Theorem \ref{pi_flat_thm} can be extended to linear systems with multiple delays. As an example, we consider the system of vibrating string with two controls proposed in \cite{Mounier_98}, which can be modeled as a set of one-dimensional wave equations together with a second order linear ordinary differential equation describing the motion of the mass. Using Mikusi\'{n}ski operational calculus (see for instance \cite{Fliess_96}), this infinite-dimensional system can be transformed into the time-delay system
\begin{equation}
\left\{\begin{array}{lll}
\psi_{1}(t)+\phi_{1}(t)-\psi_{2}(t)-\phi_{2}(t)=0,\\
\dot{\psi}_{1}(t)+\dot{\phi}_{1}(t)+\eta_{1}(\phi_{1}(t)-\psi_{1}(t))-\eta_{2}(\phi_{2}(t)-\psi_{2}(t))=0,\\
\phi_{1}(t-2\tau_{1})+\psi_{1}(t)=u_{1}(t-\tau_{1}),\\
\phi_{2}(t)+\psi_{2}(t-2\tau_{2})=u_{2}(t-\tau_{2}).\\
\end{array}\right.\label{sys_string_2control}
\end{equation}
where $\eta_{1}$ and $\eta_{2}$ are constant parameters. Denoting the state $x=(\psi_{1},\phi_{1},\psi_{2},\phi_{2})^{T}$, the control input $u=(u_{1},u_{2})$, and $\delta_{1}$, $\delta_{2}$ the delay operators of respective lengths $\tau_{1}$ and $\tau_{2}$, the system (\ref{sys_string_2control}) may be rewritten in the form $Ax=Bu$, with $A\in\mmdd{4}{4}$ and $B\in\mmdd{4}{2}$ given by
\begin{equation}
A=\left(
\begin{array}{cccc}
1&1&-1&-1\\\ddt+\eta_{1}&\ddt-\eta_{1}&\eta_{2}&-\eta_{2}\\1&\delta_{1}^{2}&0&0\\0&0&\delta_{2}^{2}&1\end{array}\right), \quad 
B=\left(\begin{array}{cc}0&0\\0&0\\\delta_{1}&0\\0&\delta_{2}\end{array}\right).
\end{equation}
The computation of a Smith-Jacobson decomposition of $B$ is here straightforward: it suffices to exchange the two last lines of $B$ with the two first lines, and we get $MBN=\left(\begin{array}{c}I_{2}\\0_{2,2}\end{array}\right)$ with
\begin{equation}
M=\left(\begin{array}{cccc}
0&0&1&0\\0&0&0&1\\\delta_{1}^{-1}&0&0&0\\0&\delta_{2}^{-1}&0&0
\end{array}\right),\quad N=I_{2}.
\end{equation}
We have then $\pi_{M,N}=\delta_{1}\delta_{2}.$
Thus $B$ is hyper-regular and
\begin{equation}
F=(0_{2,2}\quad I_{2})MA=\left(\begin{array}{cccc}\delta_{1}^{-1}&\delta_{1}^{-1}&-\delta_{1}^{-1}&-\delta_{1}^{-1}\vspace{0.3 cm}\\
\delta_{2}^{-1}\left(\ddt+\eta_{1}\right)&\delta_{2}^{-1}\left(\ddt-\eta_{1}\right)&\delta_{2}^{-1}\eta_{2}&-\delta_{2}^{-1}\eta_{2}\end{array}\right),\nonumber
\end{equation}
A right Smith-Jacobson decomposition of $F$, namely $VF\wt{Q}=\left(I_{2},\; 0_{2,2}\right)$, is given by
\begin{equation}
\begin{array}{c}
V=\left(\begin{array}{cc}\delta_{1}&0\\\delta_{2}\left(-\ddt-\eta_{1}\right)&\delta_{2}\end{array}\right)
\\
\wt{Q}=\left(\begin{array}{cccc}1&\frac{1}{2\eta_{1}}&
\frac{1}{2\eta_{1}}\left( -\ddt+ \left(\eta_{1}-\eta_{2}\right)\right)&
\frac{1}{2\eta_{1}}\left( -\ddt+ \left(\eta_{1}+\eta_{2}\right)\right)\\
0&-\frac{1}{2\eta_{1}}&
\frac{1}{2\eta_{1}}\left( \ddt+ \left(\eta_{1}+\eta_{2}\right)\right)&
\frac{1}{2\eta_{1}}\left( \ddt+ \left(\eta_{1}-\eta_{2}\right)\right)
\\0&0&1&0\\0&0&0&1\end{array}\right)
\end{array}
\end{equation}
with $\bar{\pi}=\delta_{1}\delta_{2}$, 
and where $\wt{Q}$ is obtained as the product of elementary actions $\wt{Q}_{1}$ and $\wt{Q}_{2}$:
\begin{equation}
\begin{array}{c}
\wt{Q}_{1}=\left(\begin{array}{cccc}1&-1&1&1\\0&1&0&0\\0&0&1&0\\0&0&0&1\end{array}\right),
\\
\wt{Q}_{2}=\left(\begin{array}{cccc}
1&0&0&0\\
0&-\frac{1}{2\eta_{1}}&
\frac{1}{2\eta_{1}}\left( \ddt+ \left(\eta_{1}+\eta_{2}\right)\right)&
\frac{1}{2\eta_{1}}\left( \ddt+ \left(\eta_{1}-\eta_{2}\right)\right)\\
0&0&1&0\\0&0&0&1\end{array}\right)
\end{array}
\end{equation}
thus showing that $F$ is hyper-regular.

According to step 2 of the algorithm, we compute $Q = \wt{Q}\left( \begin{array}{c}0_{2,2}\\I_{2}\end{array}\right) $ and $R = N\left( I_{2},\; 0_{2,2}\right) MAQ$:
\begin{equation}\label{QandR-ex3}
\begin{array}{c}
Q =
\left(\begin{array}{cc}
\frac{1}{2\eta_{1}}\left( -\ddt+ \left(\eta_{1}-\eta_{2}\right)\right)
&
\frac{1}{2\eta_{1}}\left( -\ddt+ \left(\eta_{1}+\eta_{2}\right)\right)
\\
\frac{1}{2\eta_{1}}\left( \ddt+ \left(\eta_{1}+\eta_{2}\right)\right)
&
\frac{1}{2\eta_{1}}\left( \ddt+ \left(\eta_{1}-\eta_{2}\right)\right)
\\
1&0
\\
0&1
\end{array}\right)
\\
R =\left(\begin{array}{cc}
R_{1,1}
&
R_{1,2}
\\
\delta_{2}^{2}
&
1
\end{array}\right)
\end{array}
\end{equation}
with
\begin{equation}
\begin{array}{l}
R_{1,1}=
\frac{1}{2\eta_{1}}\left( -\ddt+ \left(\eta_{1}-\eta_{2}\right)\right)+
\frac{\delta_{1}^{2}}{2\eta_{1}}\left( \ddt+ \left(\eta_{1}+\eta_{2}\right)\right)
\\
R_{1,2}=
\frac{1}{2\eta_{1}}\left( -\ddt+ \left(\eta_{1}+\eta_{2}\right)\right)+
\frac{\delta_{1}^{2}}{2\eta_{1}}\left( \ddt+ \left(\eta_{1}-\eta_{2}\right)\right)
\end{array}
\end{equation}
We indeed have $\pi_{R}=1$.

Therefore, setting $x=Qy$ and $u=Ry$, we obtain the expressions 
\begin{equation}
\left(\begin{array}{cc}\psi_{1}(t)\\\phi_{1}(t)\\\psi_{2}(t)\\\phi_{2}(t)\end{array}\right)
=\left(\begin{array}{cc}
\frac{1}{2\eta_{1}}\left(-\dot{y}_{1}(t)+(\eta_{1}-\eta_{2})y_{1}(t)-\dot{y}_{2}(t)+(\eta_{1}+\eta_{2})y_{2}(t)\right)
\\
\frac{1}{2\eta_{1}}\left(\dot{y}_{1}(t)+(\eta_{1}+\eta_{2})y_{1}(t)+\dot{y}_{2}(t)+(\eta_{1}-\eta_{2})y_{2}(t)\right)
\\
y_{1}(t)
\\
y_{2}(t)
\end{array}\right).
\end{equation}
and
\begin{equation}
\begin{array}{l}
u_{1}(t)=\frac{1}{2\eta_{1}}\bigl[ \dot{y}_{1}(t-\tau_{1})-\dot{y}_{1}(t+\tau_{1})+\dot{y}_{2}(t-\tau_{1})-\dot{y}_{2}(t+\tau_{1}) \bigr.
\\
\hspace{1cm}\bigl. +(\eta_{1}+\eta_{2})(y_{1}(t-\tau_{1})+y_{2}(t+\tau_{1}))+(\eta_{1}-\eta_{2})(y_{1}(t+\tau_{1})+y_{2}(t-\tau_{1}))\bigr]
\\
u_{2}(t)=y_{1}(t-\tau_{2})+y_{2}(t+\tau_{2}).
\end{array}
\end{equation}

Further, according to step 3, we compute $\wt{P}$ and $W$ of a Smith-Jacobson decomposition of $Q$, namely 
$\wt{P}QW=\left( \begin{array}{c}I_{2}\\ 0_{2,2}\end{array}\right)$.
We find $W=I_{2}$ and $\wt{P}=\wt{P}_{4}\wt{P}_{3}\wt{P}_{2}\wt{P}_{1}$ with
\begin{equation}
\wt{P}_{1}=\left(\begin{array}{cccc}0&0&1&0\\0&1&0&0\\1&0&0&0\\0&0&0&1\end{array}\right),\quad 
\wt{P}_{2}=\left(\begin{array}{cccc}
1&0&0&0\\
-\frac{1}{2\eta_{1}}\left( \ddt+\left(\eta_{1}+\eta_{2}\right)\right)&1&0&0\\
\frac{1}{2\eta_{1}}\left( \ddt-\left(\eta_{1}-\eta_{2}\right)\right)&0&1&0\\0&0&0&1\end{array}\right),\nonumber
\end{equation}
\begin{equation}
\wt{P}_{3}=\left(\begin{array}{cccc}1&0&0&0\\0&0&0&1\\0&0&1&0\\0&1&0&0\end{array}\right),\quad 
\wt{P}_{4}=\left(\begin{array}{cccc}1&0&0&0\\0&1&0&0\\
0&\frac{1}{2\eta_{1}}\left( \ddt-\left(\eta_{1}+\eta_{2}\right)\right)&1&0\\
0&-\frac{1}{2\eta_{1}}\left( \ddt+\left(\eta_{1}-\eta_{2}\right)\right)&0&1\end{array}\right),\nonumber
\end{equation}
thus
\begin{equation}
\wt{P}=\left(\begin{array}{cccc}0&0&1&0\\0&0&0&1\\
1&0&\frac{1}{2\eta_{1}}\left( \ddt-\left(\eta_{1}-\eta_{2}\right)\right)&
\frac{1}{2\eta_{1}}\left( \ddt-\left(\eta_{1}+\eta_{2}\right)\right)\\
0&1&-\frac{1}{2\eta_{1}}\left( \ddt+\left(\eta_{1}+\eta_{2}\right)\right)&
-\frac{1}{2\eta_{1}}\left( \ddt+\left(\eta_{1}-\eta_{2}\right)\right)
\end{array}\right)
\end{equation}
and 
\begin{equation}
P=W\left(I_{2}\;\; 0_{2,2}\right)\wt{P}=\left(\begin{array}{cccc}0&0&1&0\\0&0&0&1\end{array}\right)
\end{equation}

Finally, setting $y=Px$, we get $y_{1}=\psi_{2}$, $y_{2}=\phi_{2}$ and $\pi_{P}=1$.

Taking the least common multiple of $(1,1, \delta_{1}\delta_{2})$, we get $\pi=\delta_{1}\delta_{2}$. It is then immediately seen that $y_{1}=\psi_{2}$ and $y_{2}=\phi_{2}$ is a $\delta_{1}\delta_{2}$-flat output.
\begin{rem}
For multiple delays $\delta_{1}, \delta_{2},\ldots, \delta_{n}$, a flat output for which the polynomial $\pi$ is restricted to a monomial (i.e. of the form $\delta_{1}^{s_{1}}\delta_{2}^{s_{2}}\cdots \delta_{n}^{s_{n}}$) is called $\delta$-flat output in \cite{Rudolph_03}. This is the case here with $n=2$ and $s_1=s_2=1$.
\end{rem}
\begin{rem}
In \cite{Mounier_98}, a different solution $y_1=\delta_1\phi_1-u$, $y_2=\phi_1+\psi_1$ has been proposed.
\end{rem}

\section{Concluding Remarks}
In this paper, a direct characterization of $\pi$-flat outputs for linear time-varying, time-delay systems, with coefficients that are meromorphic functions of time is obtained, yielding a constructive algorithm for their  computation. The proposed approach is based on the Smith-Jacobson decomposition of a polynomial matrix over the Principal Ideal Domain $\kk(\delta)[\ddt]$, containing the original ring of multivariate polynomials $\kk[\delta, \ddt]$.
The fact that the computations are done in a larger ring, which is a principal ideal ring, makes them elementary and their localization at the powers of a $\delta$-polynomial results from an easy calculation of least common multiple of a finite set of polynomials of $\delta$. It is remarkable, however, that the $\pi$-flatness criterion only involves properties of the system matrices over the extended ring  $\kk(\delta)[\ddt]$. 
Several examples are presented to illustrate the simplicity of the approach.
Translating our algorithm in a computer algebra programme, e.g. in Maple or Mathematica, might be relatively easy and will be the subject of future works.

\section*{Acknowlegements}
The authors are indebted to Hugues Mounier and Alban Quadrat for useful discussions.

\bigskip

\appendix
\label{appendix}
\section*{Appendix}

We recall from section~\ref{alg_sec} the following theorem:
\begin{thm}[Smith-Jacobson decomposition \cite{Cohn_85,Gantmacher_66}]
Let $M\in\mmd{p}{q}$ be an arbitrary polynomial matrix of size $p\times q$. There exist unimodular matrices $U\in\uud{p}$ and $V\in\uud{q}$ such that:
\begin{itemize}
\item $\ds UMV=(\Delta_{p}|0_{p,q-p})$ if $p\leq q$,\vspace{0.5em}
\item $UMV=\left(\begin{array}{c}\Delta_{q}\\0_{p-q,q}\end{array}\right)$ if $p>q$.\vspace{0.5em}
\end{itemize}
In both cases, $\Delta_{p}\in\mmd{p}{p}$ and $\Delta_{q}\in\mmd{q}{q}$ are diagonal matrices whose diagonal elements $(d_{1},\ldots,d_{\sigma},0,\ldots,0)$ are such that $d_{i}$ is a nonzero $\ddt$-polynomial for $i=1,\ldots,\sigma$, with coefficients in $\kk(\delta)$, and is a divisor of $d_{j}$ for all $i\leq j\leq\sigma$.
\end{thm}

\section{Elementary Actions and Unimodular Matrices}

The group of unimodular matrices admits a finite set of generators corresponding to the
following \emph{elementary right and left actions}:  
\begin{itemize}
\item \emph{right actions} consist in permuting two columns, right multiplying a column by a non zero function of $\kk(\delta)$, or adding the $j$th column right multiplied by an arbitrary element of $\kk(\delta)[\ddt]$ to the $i$th column, for arbitrary $i$ and $j$; 
\item \emph{left actions} consist, symmetrically, in permuting two lines,
left multiplying a line by a non zero function of $\kk(\delta)$, or adding the
$j$th line left multiplied by an arbitrary  element of $\kk(\delta)[\ddt]$ to the $i$th line, for
arbitrary $i$ and $j$.
\end{itemize}
Every elementary action may be represented by an \emph{elementary unimodular matrix} of the form $T_{i,j}(p)=I_{\nu}+1_{i,j}p$ with $1_{i,j}$ the matrix made of a single $1$ at the intersection of line $i$ and column $j$, $1\leq i,j\leq \nu$, and zeros elsewhere, with $p$ an arbitrary  element of $\kk(\delta)[\ddt]$, and with $\nu=m$ for right actions and $\nu= n$ for left actions. 
One can easily prove that:
\begin{itemize}
\item right multiplication $MT_{i,j}(p)$ consists in adding the $i$th column  of $M$ right multiplied by $p$ to the $j$th column of $M$, the remaining part of $M$ being left unchanged, 
\item left multiplication $T_{i,j}(p)M$ consists in adding the $j$th line of $M$ left multiplied by $p$ to the $i$th line of $M$, the remaining part of $M$ being left unchanged,
\item $T_{i,j}^{-1}(p)=T_{i,j}(-p)$,
\item $T_{i,j}(1)T_{j,i}(-1)T_{i,j}(1)M$ (resp. $MT_{i,j}(1)T_{j,i}(-1)T_{i,j}(1)$) is the permutation matrix replacing the $j$th line of $M$ by the $i$th one and replacing the $j$th one of $M$ by the $i$th one multiplied by $-1$, all other lines remaining unchanged (resp. the permutation matrix replacing the $i$th column of $M$ by the $j$th one multiplied by $-1$ and replacing the $j$th one by the $i$th one, all other columns remaining unchanged).
\end{itemize}
Every unimodular matrix $V$ (left) and $U$ (right) may be obtained as
a product of such elementary unimodular matrices, possibly with a diagonal matrix 
$D(\alpha)=\diag{\alpha_{1},\ldots, \alpha_{\nu}}$ 
with $\alpha_{i}\in \kk$, $\alpha_{i}\neq 0$, $i=1,\ldots,\nu$, at the end since $T_{i,j}(p)D(\alpha)=D(\alpha)T_{i,j}(\frac{1}{\alpha_{i}}p\alpha_{j})$.

In addition, every unimodular matrix $U$ is obtained by such a product: its decomposition yields
$VU=I$ with $V$ finite product of the $T_{i,j}(p)$'s and a diagonal matrix. Thus, since the inverse of any $T_{i,j}(p)$ is of the same form, namely $T_{i,j}(-p)$, and since the inverse of a diagonal matrix is diagonal, it results that $V^{-1}=U$ is a product of elementary matrices of the same form, which proves the assertion.

\section{The Smith-Jacobson Decomposition Algorithm}\label{SJalgo-sec}

The Smith-Jacobson decomposition algorithm of the matrix $M$ consists
first in permuting columns (resp. lines) to put the element of lowest degree in upper left
position, denoted by
$m_{1,1}$, or creating this element by euclidean division (in $\kk(\delta)[\ddt]$) of two or more elements of the first line (resp. column) by suitable right actions (resp. left action). Then right divide all the other elements $m_{1,k}$  (resp. left divide the $m_{k,1}$) of the new first line (resp. first column) by $m_{1,1}$. If one of the rests is non zero, say $r_{1,k}$ (resp. $r_{k,1}$), subtract the corresponding column (resp. line) to the first column (resp. line) right multiplied (resp. left) by the corresponding quotient $q_{1,k}$ defined by the right euclidean division $m_{1,k}=m_{1,1}q_{1,k}+r_{1,k}$ (resp. $q_{k,1}$ defined by $m_{k,1}=q_{k,1}m_{1,1}+r_{k,1}$). 
Then right multiplying all the columns by the corresponding quotients
$q_{1,k}$, $k=2,\ldots, \nu$ (resp. left multiplying lines by $q_{k,1}$, $k=2,\ldots,
\mu$), we iterate this process with the transformed first line (resp. first column) until it becomes $\left(m_{1,1},0,\ldots,0\right)$ (resp.  $\left(m_{1,1},0,\ldots,0\right)^{T}$ where $^{T}$ means transposition).
We then apply the same algorithm to the second line starting from $m_{2,2}$ and so on. To each transformation
of lines and columns correspond a left or right elementary unimodular matrix and
the unimodular matrix $V$ (resp. $U$) is finally obtained as the product of
all left (resp. right) elementary unimodular matrices so constructed.

\bibliography{biblio}

\end{document}